\documentclass[12pt]{article} 

\usepackage{latexsym,amssymb, amsfonts, amscd, multicol}
\usepackage{graphics, epsfig}
\usepackage[scanall]{psfrag}
\newtheorem{definition}{Definition}

\newtheorem{fact}[definition]{Fact}

\date{}

\begin{document}

\title{Solving the Quartic with a Pencil}
\author{Dave Auckly}
\maketitle
%

%
\section{INTRODUCTION.}
\label{intro}

It is a safe bet that everyone reading this is familiar with the quadratic formula. 
\begin{fact}\label{fact1}
If $ax^2+bx+c=0$ and $a\ne 0$, then one of the following holds:
$$x=\frac{-b+\sqrt{b^2-4ac}}{2a}, $$ $$ x=\frac{-b-\sqrt{b^2-4ac}}{2a}.$$
\end{fact}
One might say that this formula allows one to solve the quadratic with a pencil. There is an analogous formula for the general quartic equation, $$ax^4+bx^3+cx^2+dx+e=0\,.$$ 
By this, we really mean four different formulas each of which gives one root of the equation. Each formula is expressible using only the operations of addition, subrtaction, multiplication, division, and extraction of roots.
If you just want to know the formula or one way to derive it--and don't care about anything else--you can just skip to the last section of this paper. The formulas for the roots of a general quartic are listed and derived there. The derivation requires the solution of the general cubic (for which we give only  hints at the derivation). Readers should be forewarned: this paper is a bit like a tourist trap. There is a main attraction, but you don't get to see it until you pass many of the souvenirs that are available along the way. The real goal of the paper is to expose readers to a number of mathematical tidbits related to the solution of the general quartic.
  
There is a mathematical notion of a ``pencil" that is rather cool and has gained prominence in geometry and topology in recent years. Pencils were also studied extensively  by algebraic geometers in the nineteenth century.
The geometric picture behind the solution to the quartic presented in the last section is a particular pencil associated with the quartic. We begin  with a description of a pencil from the view point of a topologist.

\begin{figure}
\psfrag{l0}[l][scale=.1]{\tiny $\lambda=0$}
\psfrag{ln2}[l][scale=.1]{\tiny $\lambda=-2$}
\psfrag{ln5}[l][scale=.1]{\tiny $\lambda=-5$}
\psfrag{li}[c][scale=.1]{\tiny $~\ \lambda=\infty$}
\psfrag{l5}[l][scale=.1]{\tiny $\lambda=5$}
\psfrag{l2}[l][scale=.1]{\tiny $\lambda=2$}
\psfrag{l1}[l][scale=.1]{\tiny $\lambda=1$}
\hskip70pt\epsfig{file=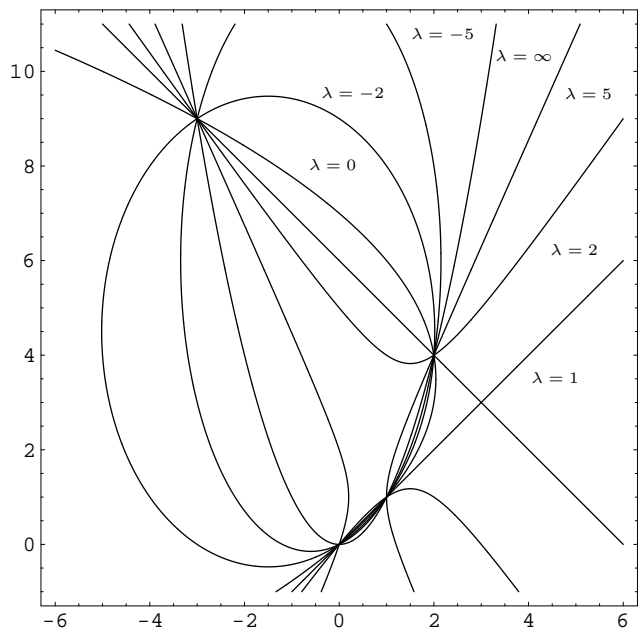,width=3.3truein}
\caption{Curves in pencil (\ref{p1}).}\label{fig1}
\end{figure}

\section{PENCILS.}
The family of quadratic sections given by
\begin{equation}\label{p1}
y^2-7y+6x+\lambda(y-x^2)=0
\end{equation}
for the different values of the parameter $\lambda$ is an example of a pencil of curves. Several members of this family (corresponding to $\lambda=-5$, $\lambda=-2$, $\lambda=0$, $\lambda=1$, $\lambda=2$, $\lambda=5$, and $\lambda=\infty$) are graphed in Figure 1, which is the key to the solution of the general quartic given in this paper. By looking at it carefully you may even be able to figure out the general formula on your own. 
If you are asking whether $\lambda$ is complex or lies in some other field, you are not an algebraic geometer. It is standard practice in algebraic geometry to consider equations as  ``functors'' that associate a geometric object to each different field or ``set of numbers.'' The figures displaying curves are the pictures for the real case, and the later figures displaying surfaces depict the complex case.

\begin{figure}
\psfrag{l0}[l][scale=.4]{\tiny $\lambda=0$}
\psfrag{li}[l][scale=.4]{\tiny $\lambda=\infty$}
\psfrag{ln1}[c][scale=.4]{\tiny $~\ \ \ \lambda=-1$}
\psfrag{l2}[l][scale=.4]{\tiny $\lambda=2$}
\psfrag{l10}[c][scale=.4]{\tiny $~\, \lambda=10$}
\epsfig{file=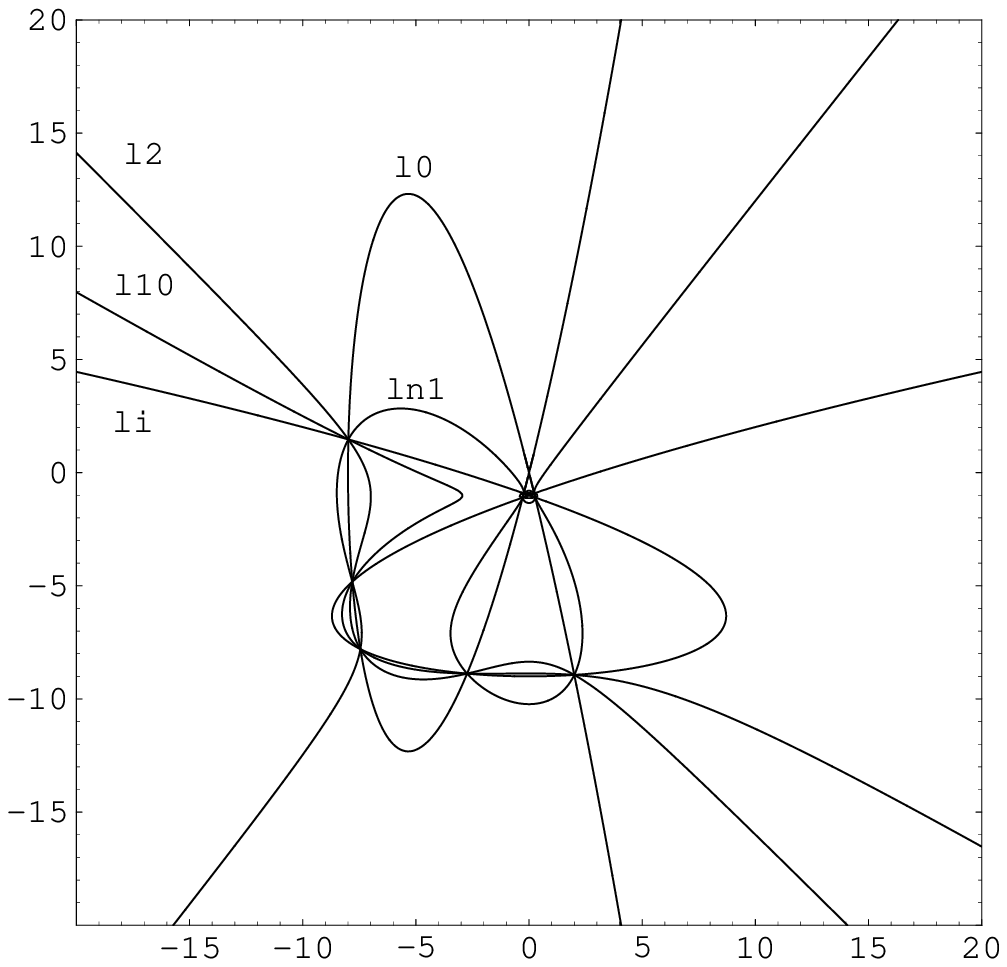,width=2.5truein}\hskip15bp\epsfig{file=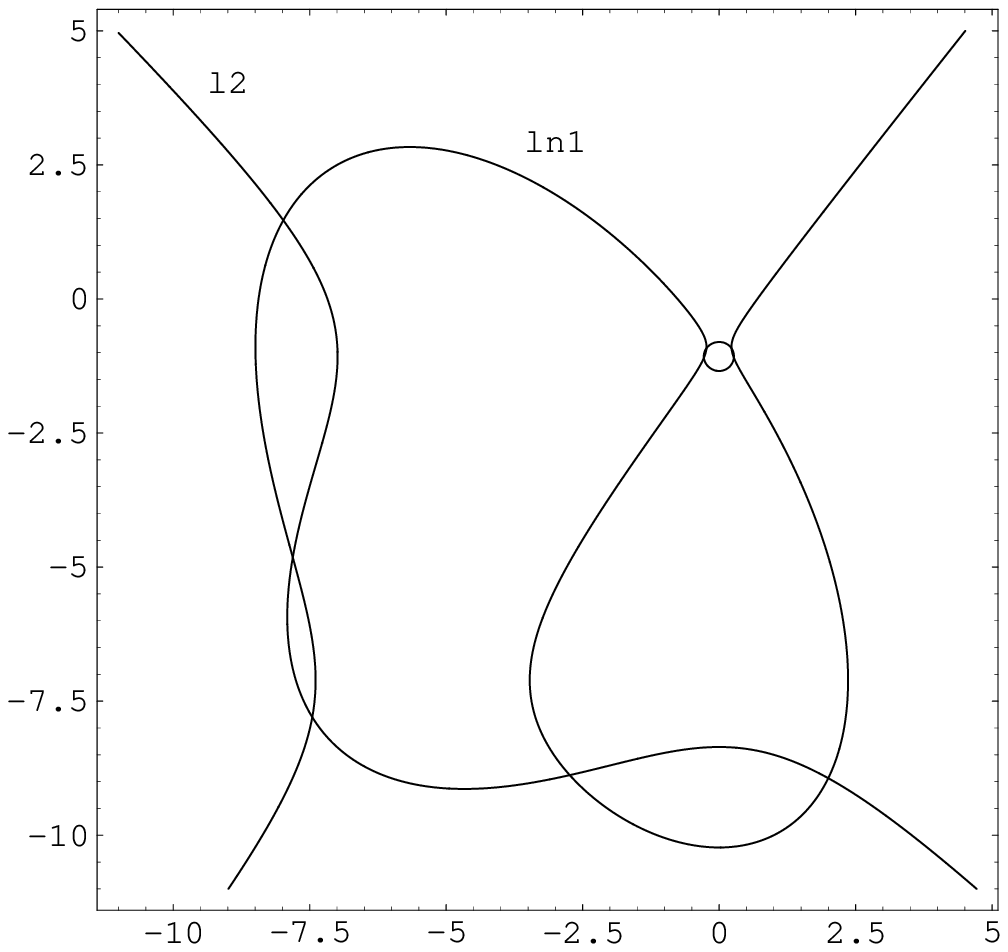,width=2.5truein}
\caption{The pencil of cubics.}\label{fig2}
\end{figure}
\begin{figure}
\psfrag{l0}[l][scale=.4]{\tiny $\lambda=0$}
\psfrag{li}[l][scale=.4]{\tiny $\lambda=\infty$}
\psfrag{ln1}[c][scale=.4]{\tiny $~\ \lambda=-1$}
\psfrag{l2}[l][scale=.4]{\tiny $\lambda=2$}
\psfrag{l10}[l][scale=.4]{\tiny $\lambda=10$}
\epsfig{file=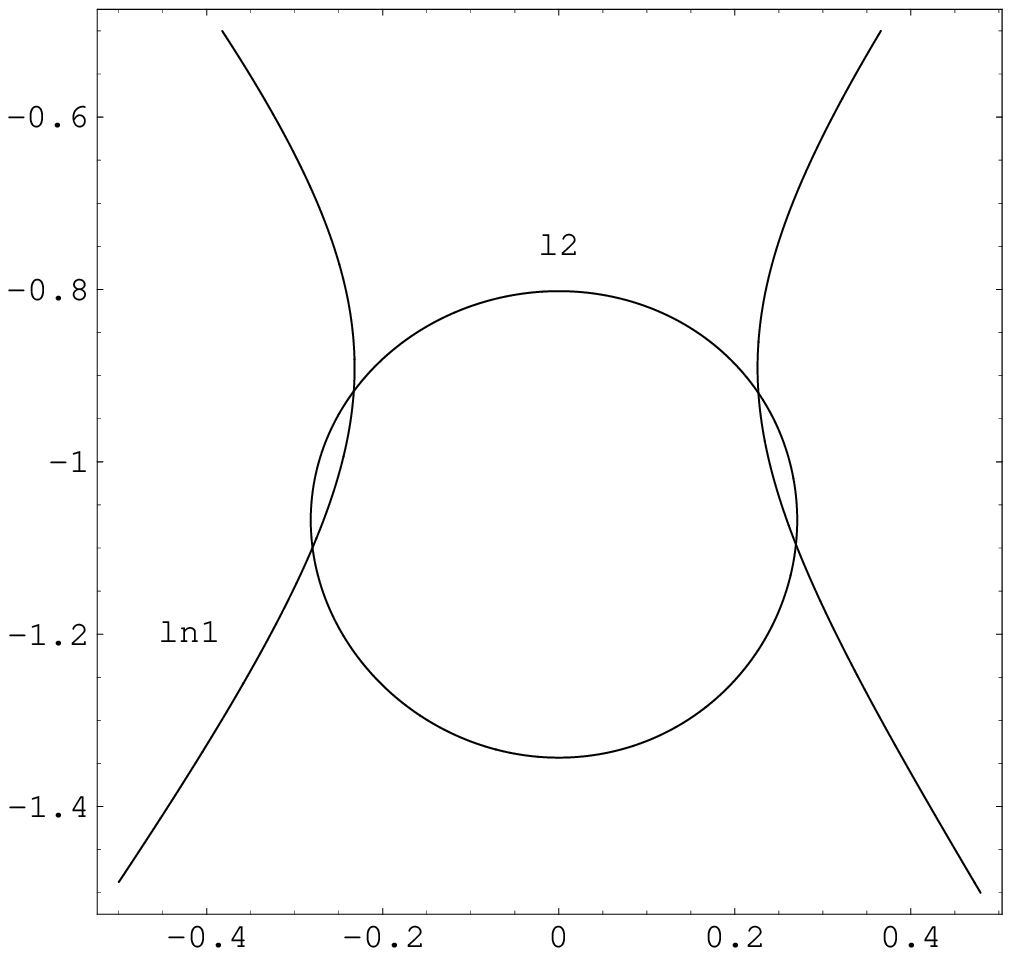,width=2.5truein}\hskip15bp\epsfig{file=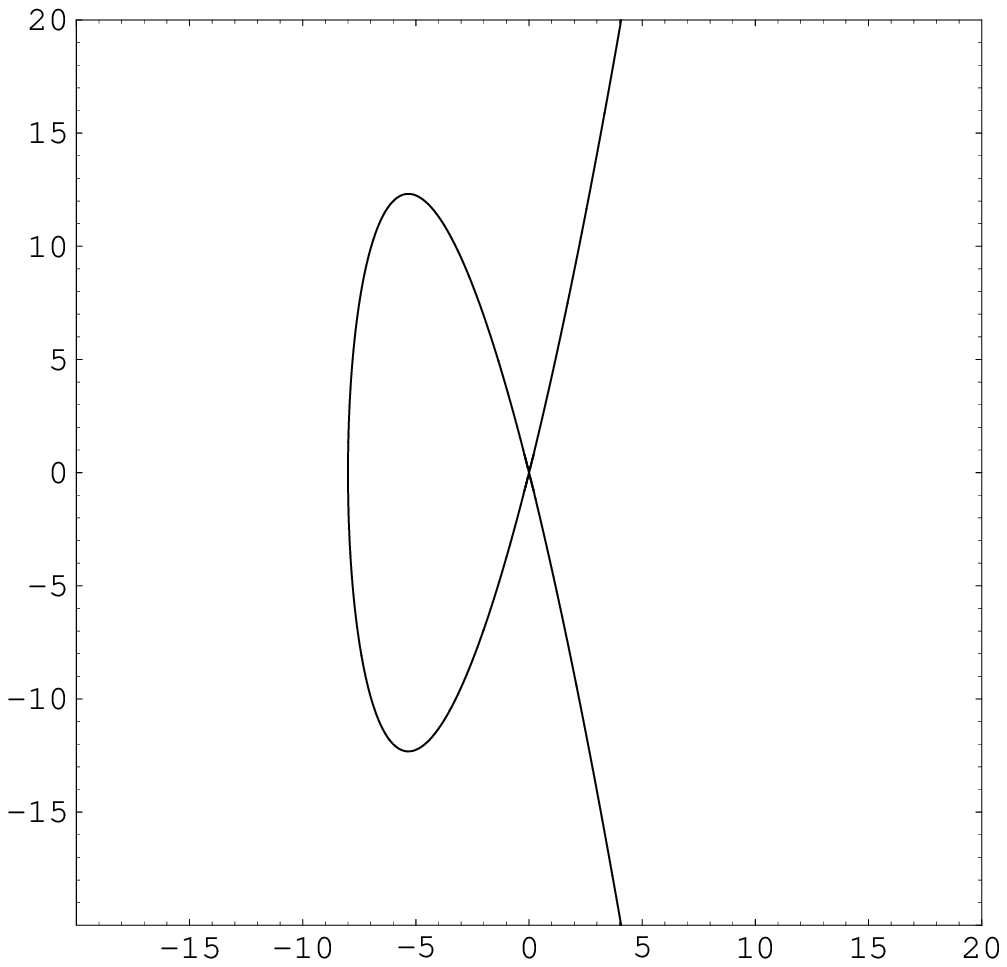,width=2.5truein}%
\caption{Curves from pencil (\ref{p2}).} 
\end{figure}

The family of curves obtained as the zeros of the family of polynomials of degree three in two variables given by 
\begin{equation}\label{p2}
y^2-2x^2(x+8)+\lambda(x^2-(y+1)^2(y+9)) =0
\end{equation}
for varying $\lambda$ is a very standard example of a pencil of curves. Several of these curves (corresponding to $\lambda=-1$, $\lambda=0$, $\lambda=2$, $\lambda=10$, and $\lambda=\infty$) are graphed in the left frame of Figure 2. Two of the generic curves (the curves with $\lambda=2$ and $\lambda=-1$) are singled out for special attention in the magnified right frame of Figure 2. If you look carefully you will see that these curves have exactly nine points in common (see the left frame of Figure 3 for a further magnification of the crossings near the origin.) Indeed, every curve in the pencil passes through each one of these nine points, which are called the {\it base points} of the pencil. For any other  point in the plane there is a unique member of the pencil that passes through it.
One of the singular curves ($\lambda=0$) is displayed in the right frame of Figure 3. The singular curve is called a {\it fish tail} by topologists and a {\it nodal cubic curve} by algebraic geometers (A formal definition of generic and singular curves is given later in this section.)
Algebraic geometers who tend to be rather smart, draw pictures of pencils similar to Figure 2, but realize that much more is happening than the picture shows directly.  

\begin{figure}
\hskip99bp\epsfig{file=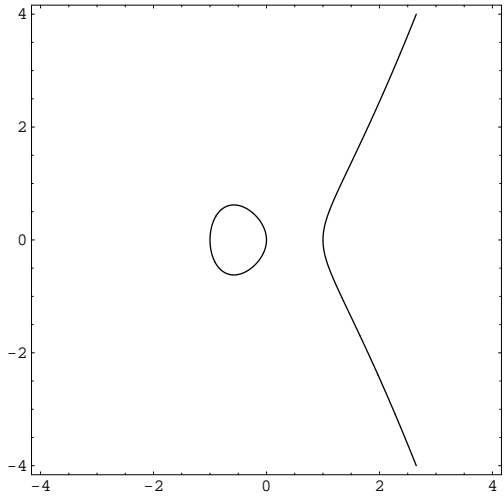,width=2.6truein}%
\caption{The curve with equation $y^2=x^3-x$.}\label{fig3}
\end{figure}

\begin{figure}
\hskip25bp\epsfig{file=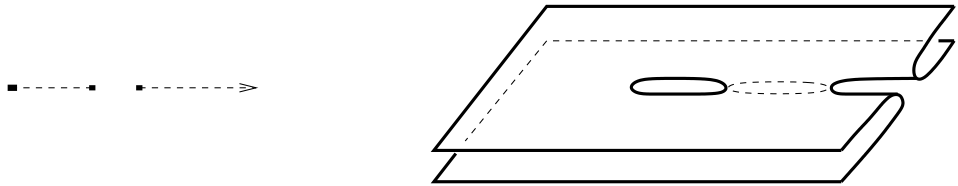,width=5truein}%
\caption{The Riemann surface of $y=\sqrt{x^3-x}$.}\label{fig4}
\end{figure}

We will now look at the geometry of pencils in greater detail, because the geometry is absolutely beautiful and really important.
The first point that arises is that it is  appropriate to use algebraically closed fields of numbers such as the field $\mathbb C$ of complex numbers. The reason for this is that  any  polynomial $p(x)=a_0+\dots +a_dx^d$ of positive degree $d$ with complex coefficients has exactly 
 $d$ zeros in ${\mathbb C}$, provided that the zeros are counted with multiplicity. The same statement cannot be made about polynomials over the real numbers. To  understand a polynomial fully one must realize that its graph is what normal people would call a surface (What an algebraic geometer calls a surface is actually four dimensional.) For example, the picture that an algebraic geometer would draw of the curve with equation $y^2=x^3-x$ is displayed in Figure \ref{fig3}. 
We can describe the actual shape of this curve by starting with the observation that, with the exception of the points $x=-1$, $x=0$, and $x=1$, every $x$-value
corresponds to two $y$-values. This means that we can construct the curve from two copies of the complex plane. The left side of Figure \ref{fig4} displays the points $x=-1$, $x=0$, and $x=1$ in the complex plane, together with a segment connecting points $x=-1$ and $x=0$ and a ray from $x=1$ to infinity. The segment and ray are called {\it branch cuts}. To go further we  use Euler's formula for the complex exponential: $e^{i\theta}=\cos\theta+i\sin\theta$. Consider the circle described by $x= e^{i(t-\pi)}/4$ for $t$ in $[0, 2\pi]$ and think about the corresponding $y$-values. We have $y=\pm\sqrt{x}\sqrt{x^2-1}$. For $t=0$ we pick $-i/2=e^{-\pi i/2}/2$ as the top sheet 
(principal branch) of $\sqrt{x}$. The values of $\sqrt{x}$ are given by $e^{{i(t-\pi)}/{2}}/2$ as $x$ traverses the circle. Notice that, after traversing the circle (i.e., after $t$ goes from $0$ to $2\pi$), $\sqrt{x}$ changes to $i/2$, which is on the other sheet (branch) (This is the first instance that we will encounter of what is called monodromy.) The value of $y$ also switches from one sheet to another as $x$ traverses the circle. We would like to match the $y$-values up consistently
with the $x$-values. This is why we introduce the branch cuts. The right side of Figure \ref{fig4} displays the result of cutting open two copies of $\mathbb C$ along the branch cuts and glueing them together in a way consistent with the equation $y^2=x^3-x$. The top sheet (principal branch) is left alone, whereas the bottom sheet (other branch) is reflected in the real axis before being glued to the principal branch. The resulting object is known as the {\it Riemann surface} for the function $y=\sqrt{x^3-x}$.
{\it Talk about this with the mathematically adept or read all about it under the heading ``Riemann surfaces'' in your library.}

A second issue that comes up is that it is useful to add points at infinity. The reason for extra points at infinity is to ensure that the zero locus of a  polynomial of degree $d$ in two variables always intersects the zero locus of a  polynomial of degree $e$ in two variables in $de$ points. Without extra points most pairs of lines in a plane intersect in one point. Properly adding extra points guarantees that {\it every} pair of distinct lines intersect in exactly one point. This is the motivation for the definition of complex projective space. 

The $n$-{\it dimensional complex projective space} is  the collection of equivalence classes of nonzero ordered $(n+1)$-tuples $(z_0,\dots,z_n)$, where two $(n+1)$-tuples are equivalent if they are complex scalar multiples of each other. (This is very similar to the definition of a rational number: the fractions $2/3$ and $4/6$ are the same.) The $n$-dimensional complex projective space is denoted by ${\mathbb C}P^n$, and its elements are denoted by $[z_0:z_1:\dots :z_n]$. The set of points of the form $[1:z_1:\dots :z_n]$ is a copy of $n$-dimensional complex space ${\mathbb C}^n$ contained in ${\mathbb C}P^n$. The {\it points at infinity} are the points of the form $[0:z_1:\dots :z_n]$. In the case of ${\mathbb C}P^1$, it is natural to think of $[w:z]$ as the fraction $z/w$, so that $[0:1]$ represents infinity. The collection of complex numbers, together with this one extra point, is known as the {\it Riemann sphere} (or {\it extended complex plane}). It is not too difficult
to extend this to a method to add points at infinity to many spaces. 

\begin{figure}
\hskip55bp\epsfig{file=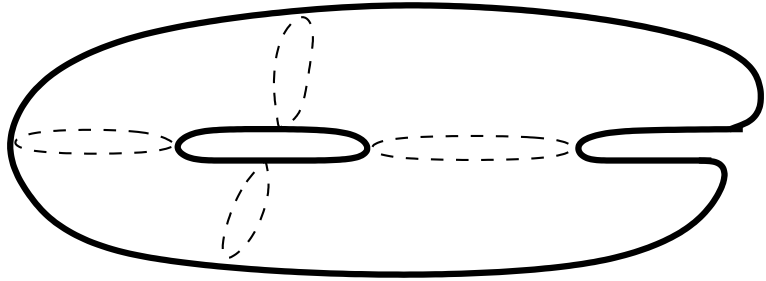,width=4truein}%
\caption{The elliptic curve with equation $ wy^2-x^3+w^2x=0$.}\label{fig5}
\end{figure}

The spaces that algebraic geometers study are the solution loci to systems of polynomial equations. If $f(z_1, \dots, z_n)$ is a polynomial function defined on ${\mathbb C}^n$, one  defines the associated homogenous polynomial by $z_0^{\hbox{deg}(f)}f(\frac{z_1}{z_0}, \dots, \frac{z_n}{z_0})$. (A function is {\it homogenous of degree} $d$ if $f(\lambda {\bf z})=\lambda^df({\bf z})$.) Points at infinity are added to the solutions of a system of polynomial equations by considering all solutions to the homogenous polynomials in ${\mathbb C}^n$. To see how this works, consider the example from Figure \ref{fig3}. The homogenous polynomial in question here is $$
w^3\left((\frac{y}{w})^2-(\frac{x}{w})^3+\frac{x}{w}\right)=wy^2-x^3+w^2x\,.
$$ 
The only zeros $[w:x:y]$ of this polynomial in ${\mathbb C}^n$ with $w=0$ have $x=0$ as well, hence there is only one such zero, namely, $[0:0:1]$.  Any zero with $w\neq 0$ may be scaled to get a zero with $w=1$. Such points correspond to the zeros of $y^2-x^3+x$. The result of glueing the one unique extra point at infinity ($[0:0:1]$) to Figure \ref{fig4} is displayed in Figure \ref{fig5}.
The process that we just described is called {\it projectivization}. {\it Talk to the math adept or look in the library for more information.}

We can now apply the projectivization proceedure to the pencil from Figure 2. Following the algebraic geometers, we will use the bizzar looking notation $H^0({\mathbb C}P^2;{\mathcal O}(L_3))$ to signify the vector space of all  homogenous polynomials of degree $3$ in three variables. Don't be intimidated by the notation. It translates to
 $$
 H^0({\mathbb C}P^2;{\mathcal O}(L_3))=\{f(x,y,z)\in {\mathbb C}[x,y,z]:f(\lambda x, \lambda y, \lambda z)=\lambda^3f(x,y,z)\}.
 $$
This is a ten-dimensional complex vector space.  (Can you find a basis for it?) 
A {\it linear system} on ${\mathbb C}P^n$ is by definition a subspace of $H^0({\mathbb C}P^n; {\mathcal O}(L))$. When one speaks of the {\it dimension} of a linear system, one usually means the dimension of the projective space associated with the linear system. This is one less than the vector space dimension of the linear system. 

\begin{definition}
A pencil on ${\mathbb C}P^n$ is a one-dimensional linear system (hence a two-dimensional linear subspace of $H^0({\mathbb C}P^n; {\mathcal O}(L))$). 
\end{definition}

Wow, that was a lot of terminology and notation to describe something fairly simple: namely, a family of polynomials of the form ``parameter times $p$ plus second parameter times $q$," in which $p$ and $q$ are linearly independent homogenous polynomials of the same degree. This entire mess exists because it helps to describe interesting features of projective algebraic varieties, and even more general spaces. The pencil in Figure 1 is in disguise. Written as a pencil, it takes the form $\mu(y^2 -7yz+6xz)+\lambda(yz-x^2)=0$. This morphs into expression (\ref{p1}) after substituting $\mu=1$ and $z=1$. The base points are the common zeros of the polynomials in the pencil.

Pencils have been defined in more general settings including spaces called projective algebraic varieties and symplectic manifolds. A {\it projective algebraic variety} is the locus in projective space of zeros to a set of homogenous polynomials.
To motivate the definition of a symplectic manifold, recall how to write Newton's equation of motion in Hamiltonian form.
Let $H$ denote the total energy of a conservative mechanical system expressed as a function of position coordinates $q_1,\dots, q_n$ and momentum coordinates $p_1,\dots, p_n$. A good example is the falling ball, which has $H=p^2/(2m)+mgq$. The equations of motion can be written in the form 
$$(\dot q,\dot p) =(\frac{\partial H}{\partial p}, -\frac{\partial H}{\partial q})=J\nabla H\,,
$$
Where $J$ is the matrix given by
$$
J=\left[\begin{array}{cc} 0 &1\\-1&0\end{array}\right]\,.
$$
Note that the first equation recovers the definition of the momentum and the second equation recovers the standard $F=ma$.
The quadratic form with matrix $J$ (generalized by the analogous block matrix for higher dimensions) is the standard {\it symplectic form}. A space obtained by identifying open subsets of ${\mathbb R}^{2n}$ via maps that preserve the standard symplectic form is called a {\it symplectic manifold}. Notice that many more maps preserve the symplectic form than the maps that preserve the dot product. A map from ${\mathbb R}^2$ to itself preserves the symplectic form if and only if it is area preserving. 
All smooth projective algebraic varieties are symplectic manifolds \cite{GH}. More information about symplectic manifolds can be found in \cite{A}. The vector space on a general
symplectic manifold $X$ corresponding to the space of homogenous polynomials of a fixed degree that one considers when
studying ${\mathbb C}P^n$ is denoted by $H^0(X; {\mathcal O}(L))$.

Returning to our example, note that
the pencil in Figure 2 is nothing more than the two-dimensional subspace of $H^0({\mathbb C}P^2;{\mathcal O}(L_3))$ given by
$$
\{ \mu(zy^2-2x^2(x+8z))+\lambda(zx^2-(y+z)^2(y+9z)):\mu, \lambda\in {\mathbb C} \}\,.
$$ 
Each curve in the figure is the locus of zeros of one of the polynomials for which $(\mu,\lambda)\ne (0,0)$. Two polynomials correspond to the same curve if and only if they are scalar multiples of each other.
 
 It is natural to label the curves associated with a pencil by elements $[\mu:\lambda]$ of one-dimensional complex projective space.
We can now start to describe the geometry associated with a pencil. It is worth considering pairs consisting of a label of a curve in the pencil together with a point on that curve. In the case of Figure 2 we obtain $$
\begin{array}{rl}
E(1)=&\!\!\!\{([\mu:\lambda],[x:y:z])\in {\mathbb C}P^1\times {\mathbb C}P^2 :\\ & \mu(zy^2-2x^2(x+8z))+\lambda(zx^2-(y+z)^2(y+9z))=0 \}\,, 
\end{array}
$$
which is a subset of ${\mathbb C}P^1\times{\mathbb C}P^2$ called the {\it rational elliptic surface}. There is a projection map $\pi_2:E(1)\to {\mathbb C}P^2$ given by $\pi_2([\mu:\lambda],[x:y:z])=[x:y:z]$. With the exception of the nine points where $zy^2-2x^2(x+8z)=0$ and $zx^2-(y+z)^2(y+9z)=0$ simultaneously, every point in ${\mathbb C}P^2$ has exactly one preimage in $E(1)$. The inverse image of any of the nine intersection points is a copy of the Riemann sphere ${\mathbb C}P^1$. We say that $E(1)$ is the result of ``blowing up" ${\mathbb C}P^2$ at nine points, while ${\mathbb C}P^2$ is the result of ``blowing down" nine spheres in $E(1)$. Algebraic geometers say that $E(1)$ and ${\mathbb C}P^2$ are {\it birationally equivalent} because they are related by a finite number of blow-ups. 

To get deeper into the geometry, consider the other projection map $\pi_1:E(1)\to {\mathbb C}P^1$, which is given by $\pi_1([\mu:\lambda],[x:y:z])=[\mu :\lambda]$. The critical points of $\pi_1$ are the points where its derivative is equal to zero. (In general {\it critical points} are the points where the derivative of a mapping is not surjective.) To  interpret this in the present context, recall that $[x:y:z]$ and $[\mu :\lambda]$ are equivalence classes. Provided that $z\ne 0$ and $\mu\ne 0$ we may represent these points in the form $[x:y:1]$ and $[1:\lambda]$. Substituting
$z=1$ and $\mu=1$ into the equations that define $E(1)$ yields 
$$\lambda=\frac{y^2-2x^2(x+8)}{(y+1)^2(y+9)-x^2}.
$$
The critical points with $z\ne 0$ and $\mu\ne 0$ can now be found by setting the partial derivatives of $\lambda$ with respect to $x$ and $y$ equal to zero. Critical points with $z=0$ or $\mu=0$ can be found similarly. The critical points of $\pi_1$ with their corresponding $\lambda$ and $\mu$ values are listed in Table 1 (first two significant digits only).

\begin{table}\label{tbl1}
\caption{Critical points.}\label{tab1}
\centering
\begin{tabular}{|l|l|l|l|l|}
\hline
$x$ & $y$ & $z$ & $\lambda$ & $\mu$ \\
\hline
$0$ & $0$ & $1$ & $0$ & $1$ \\
\hline
$0$ & $4.8$ & $1$ & $0.05$ & $1$ \\
\hline
$0$ & $-3.8$ & $1$ & $0.35$ & $1$ \\
\hline
$-8.2$ & $-6.2$ & $1$ & $7.8$ & $1$ \\
\hline
$1.0$ & $-1.0$ & $1$ & $17$ & $1$ \\
\hline
$4.1+11i$  & $-6.3-0.016i$ & $1$ & $20+11i$ & $1$ \\
\hline
$4.1-11i$ & $-6.3+0.016i$ & $1$ & $20-11i$ & $1$ \\
\hline
$-0.5-0.87i$ & $-1-0.00046i$ & $1$ & $15-.87i$ & $1$ \\
\hline
$-0.5+0.87i$ & $-1+0.00046i$ & $1$ & $15+.87i$ & $1$ \\
\hline
$-16+0.011i$ & $3.6-110i$ & $1$ & $0.0014+0.011i$ & $1$ \\
\hline
$-16-0.011i$ & $3.6+110i$ & $1$ & $0.0014-0.011i$ & $1$ \\
\hline
$0$ & $-1$ & $1$ & $1$ & $0$ \\
\hline
\end{tabular}
\end{table}

Away from the critical values (a {\it critical value} is the image of a critical point--any other point in the codomain is a {\it regular value}) the inverse image $\pi_1^{-1}([\mu:\lambda])$ is a smooth cubic similar to the curve in Figure 6 (The inverse image of a regular value is called a {\it generic curve}.) 
The inverse image of a critical value is a fish tail (The inverse image of a critical value is called a {\it singular curve}.) For example $\pi_1^{-1}([1:0])$ is the fish tail depicted in the right frame of Figure 3. The corresponding critical point has $x=0$, $y=0$, and $z=1$ (the coordinates of the point where the tail meets the fish). The inverse image of a circle that does not enclose a critical value is just a ``cylinder,'' as in the left frame of Figure 7, where each fiber (vertical line) is really the donut shape of a smooth cubic (Figure 6). The inverse image of a circle that surrounds a critical point is a fiber bundle (i.e., ``twisted cylinder'' as in the right frame of Figure 7.) We say that the cylinder or M\"obius band {\it fibers over} the circle. A space that fibers over a surface with only nondegenerate critical points and only one critical point per critical value is called a {\it Lefschetz fibration}. (A critical point is called {\it nondegenerate} if there are coordinate systems around the point and the image of the point for which the projection map takes the form $\pi(z_1,\dots,z_n)=z_1^2+\dots+z_n^2$.) A space that can be turned into a Lefshetz fibration by blowing up a number of singularities is called a {\it Lefschetz pencil}.
\begin{figure}
\hskip50bp\epsfig{file=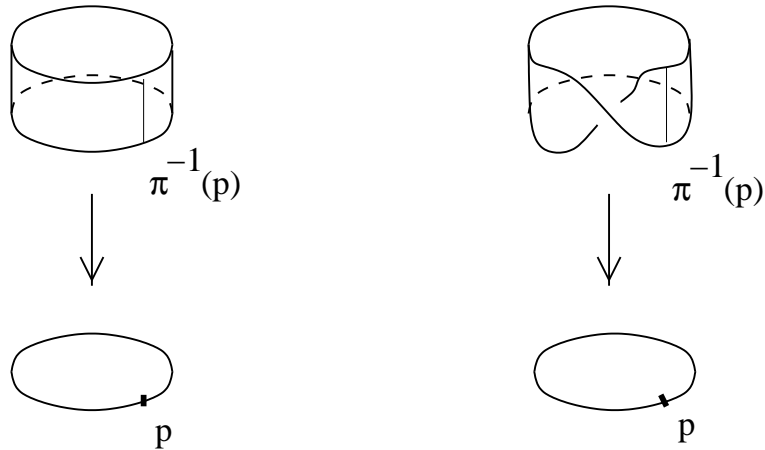,width=4truein}%
\caption{Trivial and nontrivial fibrations.}\label{fig6}
\end{figure}
Lefschetz fibrations may be described up to topological equivalence by combinatorial data that describes the monodromy (i.e., ``twisting'' around each singular fiber.) It turns out that any projective algebraic variety has the structure of (i.e., is the domain of)  a Lefschetz pencil obtained by intersecting the variety with every complex codimension one hyperplane that contains a given generic (one must avoid a small set of bad ones) complex codimension two hyperplane. The recent activity related to Lefschetz pencils stems in part from a pair of theorems, one by Simon Donaldson and the other by Bob Gompf. The theorem of Donaldson states that any symplectic manifold (an important class of objects in modern geometry) admits the structure of a Lefschetz pencil \cite{D}. Gompf's theorem establishes that any Lefschetz fibration is a symplectic manifold \cite{G}.
{\it You guessed it. Talk to the mathematically adept or read all about it at your library.} We now move from pencil propaganda to the solution of the general quartic.

\section{THE QUARTIC FORMULA.}
By scaling, any quartic equation can be written in the form
$$
x^4+ax^3+bx^2+cx+d=0.
$$
We substitute $x=u+k$ into this equation, combine like terms, and solve for the value of $k$ that forces the coefficient of $u^3$ to zero. This gives $k=-{a}/{4}$, whence $u^4+pu^2+qu+r=0$ with 
$$p=b-\frac{3a^2}{8}\,,\quad q=c-\frac{ab}{2}+\frac{a^3}{8}\,,\quad r=d-\frac{ac}{4}+\frac{a^2b}{16}-\frac{3a^4}{256}\,.
$$

Now set $v=u^2$.  Solving the system of equations
$$
v=u^2\,,\quad v^2 +pv+qu+r=0
$$ 
is easily seen to be equivalent to solving the quartic $u^4+pu^2+qu+r=0$. The quartic corresponding to the pencil in Figure \ref{fig1} is $x^4-7x^2+6x=0$. In fact, for any fixed parameter $\lambda$ the quartic $u^4+pu^2+qu+r=0$ is equivalent to the system
\begin{equation}\label{e3}
v=u^2, \qquad v^2+(p+\lambda)v-\lambda u^2 +qu+r=0.
\end{equation}
The pencil associated with Figure 1 has the form of the second equation in (\ref{e3}). 
Geometrically, the second equation represents a pencil of conic sections that share four base points. These four points provide the desired solutions to the quartic. In general, finding the points of intersection of a pair of conic sections (say, a parabola and an ellipse) is difficult. However, the pencil of the second equation in (\ref{e3}) has a singular fiber. One (corresponding to $\lambda=1$ in the specific pencil (\ref{p1})) is displayed in Figure \ref{fig1}. This singular fiber may be described as the union of two lines. Finding the points of intersection of a line and parabola is routine. There are actually two more singular fibers in this pencil. (Can you find them?)

We can rewrite the second equation in (\ref{e3}) as
$$
(v+\frac12(p+\lambda))^2-\lambda(u-\frac{q}{2\lambda})^2 = \frac14(p+\lambda)^2-\frac{q^2}{4\lambda}-r.
$$
Provided that
$$\frac14(p+\lambda)^2-\frac{q^2}{4\lambda}-r=0
$$ 
(i.e., $\lambda^3+2p\lambda^2+(p^2-4r)\lambda-q^2=0$), we can conclude that $$
(v+\frac12(p+\lambda))\pm\sqrt{\lambda}(u-\frac{q}{2\lambda})=0\,.
$$ 
Substituting $v=u^2$ into these two equations produces two quadratic equations that we can solve for $u$. The quartic formula then reads as follows:
\begin{fact}\label{fact2}
If $x^4+ax^3+bx^2+cx+d=0$, then
$$
x=\left(-\sqrt{\lambda}\pm\sqrt{\lambda-2(p+\lambda+\frac{q}{\sqrt{\lambda}})}\right)/2-{a}/{4} $$
or
$$
 x=\left(\sqrt{\lambda}\pm\sqrt{\lambda-2(p+\lambda-\frac{q}{\sqrt{\lambda}})}\right)/2-{a}/{4}\,,
$$
where
$$p=b-\frac{3a^2}{8}, \quad q=c-\frac{ab}{2}+\frac{a^3}{8}, \quad r=d-\frac{ac}{4}+\frac{a^2b}{16}-\frac{3a^4}{256},$$
and $\lambda$ is an arbitrary solution to
\begin{equation}\label{e4}
\lambda^3+2p\lambda^2+(p^2-4r)\lambda-q^2=0.
\end{equation}
\end{fact}

The  solution to the quartic equation is only useful if one has a method for solving equation (\ref{e4}). One solution to the general cubic equation is given in the following fact:
\begin{fact}\label{fact3}
If $a$, $b$, and $c$ are arbitrary complex numbers and if
\begin{equation}\label{e5}
x= {}^3\!\!\!\sqrt{-q/2+\sqrt{q^2/4+p^3/27}} +{}^3\!\!\!\sqrt{-q/2-\sqrt{q^2/4+p^3/27}} -a/3\,, 
\end{equation}
where  $p=(3b-a^2)/3$ and $q=(2a^3-9ab+27c)/27$, then $x^3+ax^2+bx+c=0$.
\end{fact}
The other solutions to the general cubic are expressible as linear combinations of the two cube roots in equation (\ref{e5}). One could try to derive the solution by constructing an appropriate pencil. Alternatively, the general cubic can be solved by making a change of variables and scale to put it in the ``reduced" form $4u^3-3u-k=0$, which can then be solved by comparison with the identity $4\cos^3\theta-3\cos\theta-\cos(3\theta)=0$. Euler's formula for $e^{i\theta}$ is used to interpret $\cos(\frac13\arccos k)$ for values of $k$ that do not lie between $\pm 1$.

Applying the foregoing method  to the quartic associated with Figure \ref{fig1},  one rewrites the pencil in the form
$$
(y+\frac12(\lambda-7))^2-\lambda(x-\frac{3}{\lambda})^2=\frac14(\lambda-7)^2-\frac{9}{\lambda}\,.
$$
Setting the right-hand side equal to zero  and simplifying gives $$
\lambda^3-15\lambda^2+49\lambda-36=0\,.
$$ 
One solution is $\lambda=1$. This leads to $(y-3)\pm (x-3)=0$, which can be combined with $y=x^2$ and solved to give the roots $0$, $1$, $2$, and $-3$ of $x^4-7x^2+6x=0$.

There is a long history associated to the solution of the quartic equation. The historical notes (written by Victor Katz) in the abstract algebra text on my shelf run from the Babylonians to the Chinese to the Arab Omar Khayyam, through several Europeans (Cardano, Tartaglia, Ferrari) and just keeps going \cite{F}. 

One of the truly amazing things about mathematics is that ideas that at first seem unrelated are often related in surprising ways. Heading out to solve the quartic equation, it is possible to get lost on a detour studying the geometry of pencils. There are many other forks in the road that we passed. One can ask, for instance, if there is a similar formula for the solution to the quintic equation. The answer is that there is no formula for the solution to the general quintic equation that can be expressed in terms of a finite number of radicals and field operations. There is, however, a solution of the general quintic expressed in terms of hypergeometric functions. This fork leads deep into algebraic territory and the realm of differential equations. The fork down the symplectic geometry path takes one deep into the mathematics of mechanical systems. 
Hopefully this little outing will encourage readers to spend some time getting lost in and picking up some souvenirs from the world of mathematics.

\medskip
\noindent{\bf ACKNOWLEDGMENTS.}
The author would like to thank the referees for numerous helpful stylistic comments.  
The  author was partially supported by 
  NSF grant DMS-0204651.

\hyphenation{Sprin-g-er-Ver-lag}

\section*{BIOGRAPHICAL SKETCH.} {\bf DAVE AUCKLY} earned his Ph.D. from the University of Michigan in 1991.
He enjoys mountaineering and spending time with his wife and son (Andrea and John).

\medskip\noindent
Department of Mathematics, Kansas State University,  
Manhattan, Kansas 66506, USA

\end{document}